\documentclass[12pt]{article}
\usepackage{amsmath,amssymb}

\def\bbd{\mathbb D}

\def\bbt{\mathbb T}

\newenvironment{pf}{{\bf Proof.~~}}{\hfill$\Box$}
\numberwithin{equation}{section}
\newtheorem{thm}{Theorem}[section]
\newtheorem{proposition}[thm]{Proposition}
\newtheorem{lemma}[thm]{Lemma}

\title{The Operator Valued Autoregressive Filter Problem and the Suboptimal
Nehari Problem in Two Variables}
\author{Jeffrey S. Geronimo\footnote{Both authors are supported in part by  grants from the
National Science Foundation and by a NATO CLG grant} and Hugo
J. Woerdeman}
\begin{document}\maketitle

\begin{center}
School of Mathematics\\
Georgia Institute of Technology\\
Atlanta, GA 30332-0160\\
geronimo@math.gatech.edu\\
and \\
Department of Mathematics\\
         Drexel University\\
         Philadelphia, PA 19104\\
Hugo.Woerdeman@drexel.edu
\end{center}


\medskip


\abstract Necessary and sufficient conditions are given for
the solvability of the operator valued two-variable autoregressive
filter problem. In addition, in the two variable suboptimal Nehari
problem sufficient conditions are given for when a strictly
contractive little Hankel has a strictly contractive symbol.
\endabstract

\section{Introduction}

The classical autoregressive filter problem asks for the
construction of an autoregressive filter based on a finite set of prescribed
correlation coefficients $c_0, \ldots , c_n$. There is a solution
to this problem if and only if the Hermitian Toeplitz matrix
$C=(c_{i-j})_{i,j=0}^n$ is positive definite, and in that case the
filter coefficients can be read off from the first column of
$C^{-1}$. While the above problem dates back to the 1950's other aspects of the theory of positive semidefinite Toeplitz
matrices had already been studied in detail in the early 1900's with the
works of Carath\'eodory, Fej\'er, Kolomogorov, Riesz, Schur,
Szeg\"o, and Toeplitz (see e.g. \cite{FF} for a full account).
Multivariable versions were considered about halfway through the
20th century. Several questions lead to extensive multivariable
results (e.g, \cite{HL,HL2}, \cite{DGK2,DGK3, DGK5}), while others
lead to counterexamples (\cite{CP}, \cite{Rud}, \cite{GK1},
\cite{Di}, \cite{LM}, \cite{L-APK}). The specific two variable
autoregressive filter problem was not completely solved until
recently in \cite{GW}. The authors found that in addition to an
expected positive definiteness requirement of a doubly
Toeplitz matrix  i.e. a block Toeplitz matrix
whose blocks are themselves Toeplitz matrices, a low rank condition on a submatrix is
necessary for the existence of a two-variable autoregressive filter with a
finite set of prescribed correlation coefficients. As it turns out,
this low rank condition may be reformulated as a commutativity
condition on matrices built form the correlation coefficients.
While this was indirectly present in the results in \cite{GW} (see
Theorem 2.2.1), it was not fully recognized as essential until
now. This commutativity condition allows for a generalization to
the operator case which we will present in this paper.

The autoregressive filter result yields sufficient conditions on a
partially defined doubly Toeplitz matrix to have a positive
definite completion, as follows. 
The notations
$ {\rm row} (c_k)_{k \in K}$ and $ {\rm col} (c_k)_{k\in K} $
stand for a row and column vector containing the entries $c_k, k \in K$,
respectively. Note that in the statement below
matrices appear that have rows and columns indexed by pairs of
integers.

\begin{thm}\label{intro} Let $c_{k}, k \in \Lambda :=\{-n, \ldots , n \}
\times \{ -m, \ldots , m \} \subset {\mathbb Z} \times {\mathbb Z}$ be given so that
$$ (c_{k-l})_{k,l \in \{ 0, \ldots , n \} \times \{ 0, \ldots , m \}} $$
is positive definite. Put
$$ \Phi = (c_{k-l})_{k,l \in \{ 0, \ldots , n-1 \} \times \{ 0, \ldots , m-1 \}}, $$
$$ \Phi_1 = (c_{k-l})_{k \in \{ 0, \ldots , n-1 \} \times \{ 0, \ldots , m-1 \},
l \in \{ 1, \ldots , n \} \times \{ 0, \ldots , m-1 \} }, $$
$$ \Phi_2 = (c_{k-l})_{k \in \{ 0, \ldots , n-1 \} \times \{ 0, \ldots , m-1 \},
l \in \{ 0, \ldots , n-1 \} \times \{ 1, \ldots , m \} }. $$
Suppose that $\Phi_1 \Phi^{-1} \Phi_2^* = \Phi_2^* \Phi^{-1}
\Phi_1$ and $$ c_{-n,m} = K_{n,m}\Phi^{-1}\tilde K_{n,m}, $$
where
\begin{equation}\label{K}
K_{n,m}={\rm row} (c_{k-l})_{k=(0,m-1), l \in \{
1, \ldots , n \} \times \{ 0, \ldots , m-1 \} },
\end{equation}
and
\begin{equation}\label{tK}
\tilde K_{n,m}={\rm
col} (c_{k-l}^*)_{k \in \{ 0, \ldots , n-1 \} \times \{ 0, \ldots
, m-1 \}, l= (n-1, 1) }.
\end{equation}
 Then there exist $c_k$, $k \not \in
\Lambda$, so that $(c_{k-l})_{k,l\in{\mathbb Z}\times {\mathbb
Z}}$ is positive definite (as an operator on $l^2 ({\mathbb
Z}\times {\mathbb Z}))$).
\end{thm}

Using the connection between positive definite and contractive
completion problems as it was used in the band method (see, e.g.,
\cite{DG5}, \cite{GKW1}, \cite{Wo2}) one may take the ideas that
go in to the proof of Theorem \ref{intro} and apply them to the
two-variable Nehari problem. The classical Nehari problem states
that a bounded Hankel operator $H  $ has an essentially
bounded symbol $\psi$, and in fact one can choose $\psi$ so that
$\| \psi \|_\infty = \| H \| $ (see, e.g., \cite{Peller}). 
In two or more variables the
situation is quite different. First of all, there are several
types of Hankels to consider. In two variables the most prominent
types are the so-called big Hankel and the little Hankel. In
\cite{CS2} it was shown that there exist bounded big Hankel
operators that do not have an essentially bounded symbol.
Recently, in \cite{FL} it was shown that every bounded small
Hankel operator has an essentially bounded symbol. The proof in
\cite{FL} relies on the dual formulation of the problem, due to
\cite{FS}. In general, though, one cannot find a symbol $\psi$ of
a small Hankel $h$, so that $\| h \| = \| \psi \|_\infty. $ We
will give sufficient conditions under which this equality can be
established in a suboptimal sense. To be more precise, we give
sufficient conditions under which $\| h \| < 1$ implies the
existence of a symbol $\psi$ so that $ \| \psi \|_\infty < 1$.

The paper is organized as follows. In Section 2 we treat the
autoregressive filter problem and as a corollary obtain Theorem
\ref{intro}. In Section 3 we treat the two-variable Nehari problem.

\section{Operator valued autoregressive filters}
\label{sec1}

A two-variable polynomial $p(z,w)$ is called {\it stable} if $p(z,w)$ is invertible for
$(z,w) \in \overline{\mathbb D}^2$, where $\overline{\mathbb D}$ stands for the closure
of ${\mathbb D} = \{ z \in {\mathbb C} : |z| < 1 \}$. Also, we denote 
${\mathbb T} = \{ z \in {\mathbb C} : |z|=1 \}$. The notation $B({\mathcal H}, {\mathcal K})$
stands for the Banach space of bounded  linear Hilbert space operators
acting ${\mathcal H} \to {\mathcal K}$. We abbreviate $B({\mathcal H}, {\mathcal H})$ as $B({\mathcal H})$. 

\begin{thm}\label{main}
Given are bounded linear operators $c_{ij} \in B({\mathcal H}),
(i,j) \in \Lambda := \{ -n,\ldots , n \} \times \{ -m, \ldots , m
\} \setminus \{ (n,m),(-n,m), (n,-m),(-n,-m)\}$. There exists
stable polynomials
\begin{equation}\label{pol} p(z,w) = \sum_{i \in \{ 0, \ldots , n \}
\atop j \in \{ 0, \ldots, m \} } p_{ij}
z^i w^j \in B({\mathcal H}), r(z,w) = \sum_{i \in \{ 0, \ldots , n \}
\atop j \in \{ 0, \ldots, m \} }
r_{ij} z^i w^j \in B({\mathcal
H}) \end{equation} with $p_{00} > 0$ and $r_{00}>0$ so that
\begin{equation}\label{cijs} p(z,w)^{*-1} p(z,w)^{-1} = \sum_{(i,j)
\in {\mathbb Z}^2 } c_{ij} z^i w^j =r(z,w)^{-1} r(z,w)^{*-1}, \
(z,w) \in {\mathbb T}^2 , \end{equation}
for some $c_{ij} \in B({\mathcal H}), (i,j) \not\in \Lambda , $ if and only if
\begin{itemize}
\item[(i)] $ \Phi_1 \Phi^{-1} \Phi_2^* = \Phi_2^* \Phi^{-1} \Phi_1 , $
\item[(ii)] when we put
 $$ c_{-n,m} = {\rm row} (c_{k-l})_{k=(0,m-1), l \in \{ 1, \ldots , n \} \times
\{ 0, \ldots , m-1 \} } \Phi^{-1} {\rm col} (c_{k-l}^*)_{k \in \{ 0, \ldots , n-1 \}
\times \{ 0, \ldots , m-1 \}, l= (n-1, 1) }, $$
then the matrices
$$ (c_{k-l})_{k,l \in \{0, \ldots , n \} \times \{ 0, \ldots , m \} \setminus \{ (n,m) \}}
\ {\rm and } \ (c_{k-l})_{k,l \in \{0, \ldots , n \} \times \{ 0, \ldots , m \} \setminus \{ (0,0) \}} $$
are positive definite.
\end{itemize}
Here
$$ \Phi = (c_{k-l})_{k,l \in \{ 0, \ldots , n-1 \} \times \{ 0, \ldots , m-1 \}}, $$
$$ \Phi_1 = (c_{k-l})_{k \in \{ 0, \ldots , n-1 \} \times \{ 0, \ldots , m-1 \},
l \in \{ 1, \ldots , n \} \times \{ 0, \ldots , m-1 \} }, $$
$$ \Phi_2 = (c_{k-l})_{k \in \{ 0, \ldots , n-1 \} \times \{ 0, \ldots , m-1 \},
l \in \{ 0, \ldots , n-1 \} \times \{ 1, \ldots , m \} }. $$
There is a unique choice for $c_{n,m}$ that results in $p_{n,m}=0$, namely
$$ c_{n,m} =
(c_{k-l})_{k = (n,m), l \in \{ 0,\ldots , n \} \times \{ 0 , \ldots , m \} \setminus \{ (0,0),
(n,m) \}}
[(c_{k-l})_{k,l\in \{ 0, \ldots , n \} \times \{ 0, \ldots , m \}
\setminus \{ (0,0), (n,m)\} }]^{-1} \times $$
$$ \ \ \ \ \ \ \ \ \ \ \ \ \ \ \ \times
(c_{k-l})_{k \in \{ 0,\ldots , n \} \times \{ 0 , \ldots , m \} \setminus \{ (0,0),
(n,m) \}, l= (0,0) }
. $$
\end{thm}

Notice that (i) is equivalent to the statement that $\Phi^{-1}
\Phi_1$ and $\Phi^{-1} \Phi_2^*$ commute. These operators
correspond exactly to the operators appearing in Theorem 2.2.1 in \cite{GW}. When
conditions (i) and (ii) are met, the polynomial $p$ may be
constructed by a Yule-Walker type of equation. Alternatively, the
Fourier coefficients $c_{ij}$ may be constructed by an iterative
process.

In the proof of Theorem \ref{main} we shall make use of some well-known results, including the
$3\times 3$ positive definite operator matrix completion problem
and the one-variable operator valued autoregressive filter
problem. We now recall these results.

\begin{proposition}\label{3x3}
Let
$$ \begin{pmatrix} A & B \cr  B^* & C \end{pmatrix} \ {\it and} \
\begin{pmatrix} C & D \cr D^* & E\end{pmatrix} $$
be positive definite Hilbert space operator matrices. Then there exist
operators $X$ so that
$$ M(X)=\begin{pmatrix} A & B & X \cr B^* & C & D \cr X^* & D^* & E \end{pmatrix} $$
is positive definite. E.g., one may choose $X= B D^{-1} C=: X_0$. In fact, $X_0$ is
the unique choice for $X$ so that $[M(X)^{-1}]_{13}=0$.
\end{proposition}

It is not hard to prove this result directly. The result also
appears in the literature in several places, e.g., in \cite{DG79},
\cite{ACC}, \cite[Section XVI.3]{FF90}. 

We will need operator valued generalizations of Theorem~2.1.5 in \cite{GW} (see also
Delsarte et al. \cite{DGK4}) and
Lemma~2.3.4 in \cite{GW}. 

\begin{thm}\label{stable}. Let
$$
p(z,w) = \sum_{i \in \{ 0, \ldots , n \}
\atop j \in \{ 0, \ldots, m \} } p_{ij}
z^i w^j \in B({\mathcal H}).
$$
Then $p(z,w)$ is stable if and only if $p(z,w)$ is invertible for all $z\in
\overline{\bbd}$ and $w\in\bbt$ and for all $z\in\bbt$ and $w\in\overline{\bbd}$.
\end{thm}
\begin{pf} Since $p(z,w)$ is invertible for all $z\in\overline{\bbd}$ and
  $w\in\bbt$ we can write
$$
p(z,w)^{-1}=\sum_{k=-\infty}^{\infty} g_k(z)w^k, \quad z\in\bbt , w \in \bbt,
$$ where $g_k(z)$ is analytic for $z\in\overline{\bbd}$. The second condition implies
that $g_k(z) =0$ for $k<0$. Thus $p(z,w)^{-1}$ is analytic for all
$(z,w)\in\overline{\bbd}^2$. Thus $p(z,w)$ is invertible for all
$(z,w)\in \overline{\bbd}^2$, and hence $p(z,w)$ is stable.
\end{pf}

\begin{lemma}\label{cholesky}
 Let $A$ be a positive definite $r \times r$ operator matrix with entries
   $A_{i,j}\in B({\mathcal H})$. Suppose that for some $1\le j < k \le r$ we have that
   $(A^{-1})_{kl}=0$, $l=1,\ldots , j$.
 Write $A^{-1}= L^*L$ where $L$ is a lower triangular operator matrix with positive 
   definite diagonal entries. Then $L$ satisfies $L_{kl}=0$, $l=1,\ldots , j$.
 Moreover, if $\tilde A$ is the $(r-1)\times (r-1)$ matrix obtained from $A$ by
 removing the $k$th row and column, and $\tilde L$ is the lower triangular
   factor 
 of $\tilde A^{-1}$ with positive diagonal entries, then
 \begin{equation} \label{LL1}
 L_{il}=\tilde L_{il}, \ i=1,\ldots , k-1; l=1,\ldots , j,
 \end{equation}
 and
 \begin{equation} \label{LL2}
 L_{i+1,l}=\tilde L_{il}, \ i=k,\ldots , r-1; l=1,\ldots , j.
 \end{equation}
 In other words, the first $j$ columns of $L$ and $\tilde L$ coincide after
 the $k$th row (which contains zeroes in columns $1,\ldots, j$) in $L$ has been removed.
 \end{lemma}

\begin{pf} Analog to the proof of Lemma~2.3.4 in \cite{GW}.
\end{pf}

\medskip

A polynomial $A(z) = A_0 +
\ldots + A_nz^n$ is called {\it stable} if $A(z)$ is invertible
for $z\in \overline{\mathbb D}$. We say that $B(z) = B_0 + \ldots + B_{-n}z^{-n}$ is {\it
antistable} if $B(1/\overline{z})^*$ is stable.

\begin{thm} (The one variable autoregressive filter problem) Let $A_j$,
$j = -n,\ldots , n$, be given Hilbert space operators, so that
the Toeplitz matrix $(A_{i-j})_{i,j=0}^n$ is positive definite. Let
$P_0, \ldots, P_n$ and $Q_{-n}, \ldots , Q_0$ be defined via
$$ \begin{pmatrix} A_0 & \cdots & A_{-n} \cr \vdots & & \vdots \cr
A_n & \cdots & A_0 \end{pmatrix} \begin{pmatrix} P_0 \cr \vdots \cr
P_n \end{pmatrix} = \begin{pmatrix} I \cr 0 \cr \vdots \cr 0 \end{pmatrix},
\begin{pmatrix} A_0 & \cdots & A_{-n} \cr \vdots & & \vdots \cr
A_n & \cdots & A_0 \end{pmatrix} \begin{pmatrix} Q_{-n} \cr \vdots \cr
Q_0 \end{pmatrix} = \begin{pmatrix} 0 \cr \vdots\cr 0 \cr I \end{pmatrix}.$$
Write $P_0 = BB^*$ and $Q_0=CC^*$ with $B$ and $C$ invertible, and put $R_i = P_i B^{*-1}$,
$S_i = Q_i C^{*-1}
$, $R (z) = \sum_{i=0}^n R_i z^i $, and $S(z) = \sum_{i=-n}^0
S_i z^i$. Then $R(z)$ is stable and $S(z)$ is anti-stable. Moreover,
$$ R(z)^{*-1} R(z)^{-1} = S(z)^{*-1} S(z)^{-1} = \sum_{j=-\infty}^\infty
A_j z^j, z \in {\mathbb T}, $$
for some $A_j=A_{-j}^*, j >n$. In fact, $A_j$, $|j|>n$, is given inductively via
$$ A_r^*= A_{-r} = \begin{pmatrix} A_{-1} & \cdots & A_{-n} \end{pmatrix} [(A_{i-j})_{i,j=0}^{n-1}]^{-1}
\begin{pmatrix} A_{-r+1} \cr \vdots \cr A_{-r+n} \end{pmatrix}, r \ge n+1. $$
\end{thm}

The matrix version of this result goes back to \cite{DG79a}. The operator valued
case appeared first in \cite{GKW2}. One may also consult \cite[Section III.3]{Wo2}
or \cite[Chapter XXXIV]{GGK2}.

We will need the notions of left and right stable factorizations
of operator valued trigonometric polynomials.  Let $A(z) =
\sum_{i=-n}^n A_i z^i$ be a matrix-valued trigonometric polynomial
that is positive definite on $\bbt$, i.e., $A(z) >0$ for $|z|=1$.
In particular, since the values of $A(z)$ on the unit circle are
Hermitian, we have $A_i = A_{-i}^*$, $i=0,\ldots , n$. The
positive matrix function $A(z)$ allows a {\it left stable
factorization}, that is, we may write
$$ A(z) = M(z) M(1/\overline{z})^*, z\in {\mathbb C} \setminus \{ 0 \} , $$
with $M(z)$ a stable matrix polynomial of degree $n$. In the
scalar case, this is the well-known Fej\'er-Riesz factorization
and goes back to the early 1900's. For the matrix case the result
goes back to \cite{RO} and \cite{HE}. When we require that $M(0)$
is lower triangular with positive definite diagonal entries, the stable
factorization is unique. We shall refer to this unique factor
$M(z)$ as {\it the left stable factor} of $A(z)$. Similarly, we
define right variations of the above notions. In particular, if
$N(z)$ is so that $A(z)= N(1/\overline{z})^* N(z), z \in {\mathbb
C} \setminus \{ 0 \}$, $N(z)$ is stable and $N(0)$ is lower
triangular with positive definite diagonal elements, then $N(z)$ is called
the {\it right stable factor} of $A(z)$.

{\bf Proof of Theorem \ref{main}.}
Observe that $\Phi_i \Phi^{-1}$ and $\Phi^{-1} \Phi_i$, $i=1,2$, have the following
companion type forms:
\begin{equation}\label{form1} \Phi_1 \Phi^{-1} = \begin{pmatrix} * & \cdots
    & * & * \cr I & 0 & & \cr & \ddots &\ddots & \cr & & I & 0
\end{pmatrix}, \Phi^{-1} \Phi_1 = \begin{pmatrix} 0 & & & * \cr I & \ddots& & * \cr & \ddots &0 & \vdots
\cr & & I & * \end{pmatrix} , \end{equation}
\begin{equation}\label{form2} \Phi_2^* \Phi^{-1} = (Q_{ij})_{i,j=0}^{n-1}, Q_{ij} = \begin{pmatrix} * & \cdots & * & * \cr
\delta_{ij} I & 0& & \cr & \ddots &\ddots & \cr & & \delta_{ij} I & 0
\end{pmatrix}, \end{equation}
\begin{equation}\label{form3}  \Phi^{-1} \Phi_2^* = (R_{ij})_{i,j=0}^{n-1}, R_{ij} = \begin{pmatrix} 0 & & & * \cr \delta_{ij} I & \ddots& & * \cr & \ddots &0 & \vdots
\cr & & \delta_{ij} I & * \end{pmatrix} , \end{equation}
where $\delta_{ij}=1$ when $i=j$ and $\delta_{ij}=0$ otherwise.
Consequently, if $S= (S_{ij})_{i,j=0}^{n-1}$ satisfies
\begin{equation}\label{blockT}
\Phi_1 \Phi^{-1} S = S \Phi^{-1} \Phi_1 ,
\end{equation}
then $S$ is block Toeplitz (i.e., $S_{ij} = S_{i+1,j+1}$ for all $0\le i,j \le n-2$).
Next, if $S=(S_{ij})_{i,j=0}^{n-1}$ satisfies
\begin{equation}\label{T}
\Phi_2^* \Phi^{-1} S = S \Phi^{-1} \Phi_2^* ,
\end{equation}
then each $S_{ij}$ is Toeplitz. It follows from (i) that all expressions
of the form
\begin{equation}\label{form}
S = \Psi_{i_1} \Phi^{-1} \Psi_{i_2} \Phi^{-1} \cdots \Phi^{-1} \Psi_{i_k},
\end{equation}
where $i_j \in \{1,2\}, \Psi_1 = \Phi_1$ and $\Psi_2 = \Phi_2^*$, satisfy
\eqref{blockT} and \eqref{T}. Thus all expressions $S$ in \eqref{form}
are doubly Toeplitz. In particular, $ \Phi_1 \Phi^{-1} \Phi_2^* = \Phi_2^* \Phi^{-1} \Phi_1  $
is doubly Toeplitz. Upon closer inspection we have that
\begin{equation}\label{SS}
\Phi_1 \Phi^{-1} \Phi_2^* = (\Gamma_{i-j})_{i=-1, j=0}^{n-2\ \ \ n-1},
\Gamma_{k} = (c_{k,r-s})_{r=1, s=0}^{m\ \ \ m-1},
\end{equation}
where $c_{-n,m}$ is defined by this equation to be as under (ii). Notice that due to
\eqref{SS} we have that
\begin{equation}\label{SSS}
\begin{pmatrix} c_{-1,1} & \cdots & c_{-n,1} \cr \vdots & & \vdots \cr c_{-1,m} & \cdots & c_{-nm}
\end{pmatrix} = \begin{pmatrix} I & 0 & \cdots & 0 \end{pmatrix} \Phi_1
\Phi^{-1} \Phi_2^* \begin{pmatrix} e_0 & & 0 \cr & \ddots & \cr 0 & & e_0 \end{pmatrix},
\end{equation}
with $e_0 = \begin{pmatrix} 1 & 0 & \cdots & 0 \end{pmatrix}^*$.
Due to (ii) and the $3\times 3$ positive definite matrix completion problem, we can choose
$c_{n,m} = c_{-n,-m}^*$ so that $\Gamma:=(c_{k-l})_{k,l \in \{ 0, \ldots , n \} \times
\{ 0, \ldots , m \} } > 0$. View $\Gamma = (C_{i-j})_{i,j=0}^n$ where
$C_k = (c_{k,r-s})_{r,s=0}^m$, and extend $\Gamma$ following the one variable
theory to $(C_{i-j})_{i,j=0}^\infty$, where
$$ C_r^*= C_{-r} = \begin{pmatrix} C_{-1} & \cdots & C_{-n} \end{pmatrix} [(C_{i-j})_{i,j=0}^{n-1}]^{-1}
\begin{pmatrix} C_{-r+1} \cr \vdots \cr C_{-r+n} \end{pmatrix}, r \ge n+1. $$
Equivalently, if we let
$$ \begin{pmatrix} Q_0 \cr \vdots \cr Q_n \end{pmatrix} = \Gamma^{-1} \begin{pmatrix} I \cr
0 \cr \vdots \cr 0 \end{pmatrix},
$$
and we factor $Q_0 = LL^*$ with $L$ lower triangular, and put $P_j = Q_j L^{*-1}$, $j=0,\ldots ,n$,
then $P(z) := P_0 + \ldots + z^n P_n $ is stable and
\begin{equation}\label{cp}
\sum_{j=-\infty}^\infty C_jz^j = P(z)^{*-1} P(z)^{-1} , z \in {\mathbb T}.
\end{equation}
Due to \eqref{SSS} it follows from
Lemma~\ref{cholesky} that $P_j$ is of the form
$$ P_j = \begin{pmatrix} p_{j0} & 0 \cr \begin{pmatrix} p_{j1} \cr \vdots \cr p_{jm} \end{pmatrix}
& \tilde{P}_j \end{pmatrix}, j=0,\ldots , n. $$
But then it follows that $\tilde P(z) :=\tilde  P_0 + \ldots + z^n \tilde P_n $ is stable, and that
$$ \sum_{j=-\infty}^\infty \tilde C_j z^j = \tilde P (z)^{*-1} \tilde P (z)^{-1},
z \in {\mathbb T}, $$
where $\tilde C_j$ is obtained from $C_j$ by leaving out the first row and column.

Similarly, if we let $$ \begin{pmatrix} R_{-n} \cr \vdots \cr R_0 \end{pmatrix} = \Gamma^{-1} \begin{pmatrix} 0 \cr
\vdots \cr 0 \cr I \end{pmatrix},
$$
and we factor $R_0 = UU^*$ with $U$ upper triangular, and put $S_j = R_j U^{*-1}$, $j=-n,\ldots ,0$,
then $S(z) := S_0 + \ldots + z^{-n} S_{-n} $ is anti-stable and
$$ \sum_{j=-\infty}^\infty C_jz^j = S(z)^{*-1} S(z)^{-1} , z \in {\mathbb T}. $$
Due to \eqref{SSS} it follows from Lemma \ref{cholesky} that $S_j$ is of the form
\begin{equation}\label{formsj}
S_j = \begin{pmatrix} \tilde S_j & \begin{pmatrix} \tilde p_{j,-m} \cr \vdots \cr \tilde p_{j,-1} \end{pmatrix} \cr 0
& \tilde{p}_{0j} \end{pmatrix}, j=-n,\ldots , 0.
\end{equation}
But then it follows that $\tilde S(z) :=\tilde  S_0 + \ldots + z^{-n} \tilde S_{-n} $ is anti-stable, and that
$$ \sum_{j=-\infty}^\infty \hat C_j z^j = \tilde S (z)^{*-1} \tilde S (z)^{-1},
z \in {\mathbb T}, $$
where $\hat C_j$ is obtained from $C_j$ by leaving out the last row and column.

Due to the block Toeplitzness of $C_j$, $j=-n,\ldots , n$, we have that
$\tilde C_j = \hat C_j$, $j=-n,\ldots , n $. As $\tilde S (z) $ and $\tilde P (z)$ follow the
one variable construction with $(\tilde C_{i-j})_{i,j=0}^n = (\hat C_{i-j})_{i,j=0}^n$, we have
by the one variable theory that
$$ \tilde P(z)^{*-1} \tilde P(z)^{-1} = \tilde S(z)^{*-1} \tilde S(z)^{-1}, \ z \in
{\mathbb T}, $$
and thus $\tilde C_j = \hat C_j$, $j \in {\mathbb Z}$. Thus $C_j$ is Toeplitz for all $j$.
Denote $C_j = (c_{j,r-s})_{r,s=0}^m.$
As $\sum_{j=-\infty}^\infty C_j z^j >0$, $j \in {\mathbb T}$, we have that the infinite
block Toeplitz $(C_{i-j})_{i,j = -\infty}^\infty$ is positive definite. We may regroup
this infinite block Toeplitz matrix with Toeplitz entries as
$(T_{i-j})_{i,j=0}^m$ where
$$ T_j = (c_{r-s,j})_{r,s=-\infty}^\infty. $$
Taking equality \eqref{cp}, and performing a regouping and extracting the first column from
$P_i$, one arrives at
$$ \begin{pmatrix} T_0 & \cdots & T_{-m} \cr \vdots & & \vdots \cr T_m & \cdots & T_0
\end{pmatrix} \begin{pmatrix} \Pi_0 \cr \vdots \cr \Pi_m \end{pmatrix} =
\begin{pmatrix} Q_0 \cr 0 \cr \vdots \cr o \end{pmatrix} , $$
where
$\Pi_j = (p_{r-s,j})_{r,s=-\infty}^\infty$, $Q_0 = \Pi_0^{-1}= (q_{r-s,0})_{r,s=-\infty}^\infty$,
$p_{ij}=0$ for $i<0$ or $i > n$, $j<0$ or $j>m$ and
$$ q(z) = \sum_{i=-\infty}^0 q_{i0}z^i := (\sum_{i=0}^n p_{i0}z^i)^{*-1}. $$
Note that $q(z)$ is indeed anti-analytic as $\sum_{i=0}^n p_{i0}z^i$ is stable. The one
variable theory now yields that
$$\Pi (w) := \Pi_0 + \ldots + \Pi_m w^m $$
is invertible for all $w \in \overline{\mathbb D}$. As $\Pi (w)$ is Toeplitz, its symbol
is invertible on ${\mathbb T}$, and thus $p(z,w) = \sum_{i=0}^n \sum_{j=0}^m p_{ij}z^i w^j$
is invertible for all $|w| \le 1$ and $|z|=1$. By reversing the roles of $z$ and $w$
one can prove in a similar way that $p(z,w)$ is invertible for all $|z| \le 1$ and
$|w|=1$. Combining these two statements yields by Theorem \ref{stable} that $p(z,w)$ is stable. In addition, we obtain that
$$ \Pi (w)^{*-1} \Pi (w)^{-1} $$
has Fourier coefficients $T_{-m}, \ldots , T_m$. But then it follows that
$$ p(z,w)^{*-1} p(z,w)^{-1} $$
has Fourier coefficients $c_{ij}$. Similarly, one proves that
$r(z,w):= \sum_{i\in \{ 0,\ldots , n\} \atop j \in \{ 0,\ldots , m \}}  \tilde{p}_{-i,-j}^* z^i w^j$ is stable and $
r(z,w)^{-1} r(z,w)^{*-1} $ has Fourier coefficients $c_{ij}$. This
proves one direction of the theorem.

For the converse, let $p$ and $r$ as in \eqref{pol} be stable and
suppose that \eqref{cijs} holds. Denote $f(z,w) = \sum_{(i,j)
\in{\mathbb Z}^2 } c_{ij} z^i w^j$. Write
$f(z,w)=\sum_{i=-\infty}^\infty f_i(z) w^i . $ Then $ T_k(z):=
(f_{i-j}(z))_{i,j=0}^k >0$ for all $k \in {\mathbb N}_0$ and all
$z\in {\mathbb T}$. Next, write
$$ p(z,w)= \sum_{i=0}^m p_i(z) w^i , r(z,w)=\sum_{i=0}^m
r_i(z) w^i , $$ and put $p_i(z) =r_i(z) \equiv 0$ for $i > m$. By
the inverse formula for block Toeplitz matrices  \cite{GH} we have
that for $k \ge m-1$ and $z\in {\mathbb T}$
\begin{equation}
\label{ttnine}
\begin{array}{rl} & T_{k}(z)^{-1} =
\begin{pmatrix} p_0 (z) &&\bigcirc\cr\vdots&\ddots&&\cr p_{k}(z)
&\cdots&p_0 (z)\end{pmatrix} \begin{pmatrix}\bar p_0
(1/\overline{z})^* & \cdots & p_{k}(1/\overline{z})^*\cr &\ddots
&\vdots\cr \bigcirc&& p_0(1/\overline{z})^*\end{pmatrix}- \cr
&-\begin{pmatrix} r_{k+1}(1/\overline{z})^* &&\bigcirc\cr \vdots
&\ddots&&\cr
 r_1(1/\overline{z})^* &\cdots& r_{k+1}(1/\overline{z})^*\end{pmatrix}
\begin{pmatrix} r_{k+1}(z) &\cdots&r_1(z)\cr & \ddots&\vdots \cr
\bigcirc && r_{k+1}(z)\end{pmatrix} =: E_k(z) \end{array}
\end{equation}
As was proven in \cite[Proposition 2.1.2]{GW} for the scalar case,
we have that for $k \ge m-1$, the left stable factors $M_k(z)$ and
$M_{k+1}(z)$ of $E_k(z)$ and $E_{k+1}(z)$, respectively, satisfy
\begin{equation}
\label{ttten} M_{k+1} (z) = \begin{pmatrix} p_0 (z) & 0 \cr {\rm col}
(p_l(z))_{l=1}^{k+1} & M_k(z) \end{pmatrix} .
\end{equation}
Indeed, if we define $M_{k+1}(z)$ by this equality, then writing
out the product $M_{k+1}(z) M_{k+1}(1/\bar z)^*$ and comparing it
to $ E_{k+1}(z)$, it is straightforward to see that $M_{k+1}(z)
M_{k+1}(1/\bar z)^* = E_{k+1}(z)$. Since both $p_0(z)$ and
$M_k(z)$ are stable, $M_{k+1}(z)$ is stable as well. Moreover,
since $p_0(0)>0$ and $M_k(0)$ is lower triangular with positive
diagonal entries, the same holds for $M_{k+1}(0)$. Thus
$M_{k+1}(z)$ must be the stable factor of $E_{k+1}(z)$. Let $C_k=
(c_{k,r-s})_{r,s=0}^m$ as before. Then we have that $T_m(z)=
\sum_{k=-\infty}^\infty C_k z^k =M_{m}(1/\bar z)^{*-1} M_m(z)$.
Writing $M_m(z) = P_0 +\ldots + z^n P_n$ it follows from the
one-variable result that
$$ \begin{pmatrix} C_0 & \cdots & C_{-n} \cr \vdots & & \vdots \cr
C_n & \cdots & C_0 \end{pmatrix} \begin{pmatrix} P_0 \cr \vdots
\cr P_n \end{pmatrix} = \begin{pmatrix} P_0^{*-1} \cr 0 \cr \vdots
\cr 0 \end{pmatrix}. $$ Due to the zeros in $P_1 , \ldots , P_n$ (see \eqref{ttten})
it follows from Proposition \ref{3x3} that \eqref{SSS} holds. By a
similar argument, reversing the roles of $z$ and $w$, we obtain that
$$ \begin{pmatrix}\widetilde C_0 & \cdots &\widetilde C_{-n} \cr \vdots & &
\vdots \cr\widetilde C_n & \cdots &\widetilde C_0
\end{pmatrix}\begin{pmatrix} S_{-n} \cr \vdots
\cr S_0 \end{pmatrix} = \begin{pmatrix} 0 \cr \vdots \cr 0 \cr
S_0^{*-1} \end{pmatrix}, $$ where $\widetilde
C_k=(c_{r-s,k})_{r,s=0}^n$ and $S_j$ has the form as in \eqref{formsj}. Using the zero
structure of $S_{-1}, \ldots , S_{-n}$ one obtains equality
\eqref{SSS} with $\Phi_1\Phi^{-1} \Phi_2^*$ replaced by $\Phi_2^*
\Phi^{-1} \Phi_1$. But then it follows that
\begin{equation}\label{eq}
\begin{pmatrix} I & 0 & \cdots & 0 \end{pmatrix} \Phi_1
\Phi^{-1} \Phi_2^* \begin{pmatrix} e_0 & & 0 \cr & \ddots & \cr 0 & & e_0 \end{pmatrix} =
\begin{pmatrix} I & 0 & \cdots & 0 \end{pmatrix} \Phi_2^*
\Phi^{-1} \Phi_1 \begin{pmatrix} e_0 & & 0 \cr & \ddots & \cr 0 &
& e_0 \end{pmatrix}. \end{equation} Due to
\eqref{form1}-\eqref{form3} it is easily seen that $\Phi_2^*
\Phi^{-1} \Phi_1$ and $\Phi_1 \Phi^{-1} \Phi_2^*$ have the same
block entries anywhere else, so combining this with \eqref{eq}
gives that  $\Phi_2^* \Phi^{-1} \Phi_1 = \Phi_1 \Phi^{-1}
\Phi_2^*$. This yields (i) and the equality for $c_{-n,m}$ in (ii). 
The positive definiteness of the
matrices in (ii) follows as they are restriction of the
multiplication operator with symbol $f$, which takes on positive
definite values on ${\mathbb T}^2$. \hfill $\Box$

{\bf Proof of Theorem \ref{intro}.} Follows directly from Theorem \ref{main}.
\hfill $\Box$

\section{Nehari's problem in two variables}

We start by stating a version of the operator valued one-variable Nehari result that will be
useful in our two-variable result. The operator valued Nehari result is due to Page
\cite{Page} who proved it using its connection to the commutant lifting theorem, and
independently to Adamjan, Arov and Krein \cite{AAK4} who had a matricial approach. The latter
approach is close to the one we employ here.

We let $l^2_{\mathcal H} (K)$ denote the Hilbert space of sequences
$\eta=(\eta_j)_{j\in K}$ satisfying $\| \eta \|:= \sqrt{\sum_{j\in K} \| \eta_j \|^2_{\mathcal H}}
< \infty $. We shall typically write Hankels in a Toeplitz like format by reversing
the order of the columns of our Hankel matrices. E.g., in the one-variable case our Hankels
shall typically act $l^2 (-{\mathbb N}_0 ) \to l^2 ({\mathbb N}_0 )$ as opposed to the usual
convention of acting $l^2 ({\mathbb N}_0 ) \to l^2 ({\mathbb N}_0 )$.

\begin{thm}\label{Nehari} Let $\Gamma_i \in L({\mathcal H}, {\mathcal K}), i \ge 0$, be bounded linear
Hilbert space operators so that the Hankel
\begin{equation}\label{hankel}
H:= \begin{pmatrix} & & \Gamma_1 & \Gamma_0 \cr & & \ddots & \Gamma_1 \cr & & &
\cr & & & \end{pmatrix} : l^2_{\mathcal K} ( - {\mathbb N}_0 ) \to l^2_{\mathcal H} ( {\mathbb N}_0 ),  \end{equation}
is a strict contraction. Solve for operators $\Delta_0, D_{-1}, D_{-2} , \ldots ,
B_0 , B_1, \ldots  $, satisfying the
Yule-Walker type equation
\begin{equation}\label{YW}
\begin{pmatrix} I & H \cr H^* & I \end{pmatrix} \begin{pmatrix} B \cr D \end{pmatrix}
= \begin{pmatrix} 0 \cr \Delta \end{pmatrix}, \end{equation}
where
$$ B =\begin{pmatrix} B_0 \cr B_1 \cr B_2 \cr \vdots \end{pmatrix} : {\mathcal K} \to
l^2_{\mathcal H} ({\mathbb N}_0) , D= \begin{pmatrix} \vdots \cr
D_{-2} \cr D_{-1} \cr I_{\mathcal K} \end{pmatrix}: {\mathcal K} \to
l^2_{\mathcal K} (-{\mathbb N}_0) , \Delta = \begin{pmatrix} \vdots \cr 0 \cr 0 \cr
\Delta_0 \end{pmatrix}: {\mathcal K} \to
l^2_{\mathcal K} (-{\mathbb N}_0) . $$
For $j=-1,-2, \ldots,$ put
\begin{equation}\label{gj}
\Gamma_j = - \Gamma_{j+1} D_{-1} - \Gamma_{j+2} D_{-2} - \ldots = - \sum_{k=1}^\infty
\Gamma_{j+k} D_{-k}. \end{equation}
Then $f \sim \sum_{j=-\infty}^\infty \Gamma_j z^j$ belongs to
$L^\infty_{\mathcal H} ({\mathbb T} ) $ and $\| f \|_\infty <1 $.

Alternatively, the Fourier coefficients $\Gamma_j$ of $f$ may be constructed as follows.
Solve for operators $\alpha_0, A_1, A_2, \ldots, C_0, C_{-1} , \ldots
$, satisfying the
Yule-Walker type equation
\begin{equation}\label{YW2}
\begin{pmatrix} I & H \cr H^* & I \end{pmatrix} \begin{pmatrix} A \cr C \end{pmatrix}
= \begin{pmatrix} \alpha \cr 0 \end{pmatrix}, \end{equation}
where
$$ A =\begin{pmatrix} I_{\mathcal H} \cr A_1 \cr A_2 \cr \vdots \end{pmatrix} : {\mathcal H} \to
l^2_{\mathcal H} ({\mathbb N}_0), C= \begin{pmatrix} \vdots \cr
C_{-2} \cr C_{-1} \cr C_0 \end{pmatrix} : {\mathcal H} \to
l^2_{\mathcal K} (-{\mathbb N}_0), \alpha = \begin{pmatrix} \alpha_0 \cr 0 \cr 0 \cr
\vdots \end{pmatrix} : {\mathcal H} \to
l^2_{\mathcal H} ({\mathbb N}_0) . $$
For $j=-1,-2, \ldots,$, $\Gamma^*_j$ may be calculated from,
\begin{equation}\label{gj2}
\Gamma_j^* = - \Gamma_{j+1}^* A_{1} - \Gamma_{j+2}^* A_{2} - \ldots = - \sum_{k=1}^\infty
\Gamma_{j+k}^* A_{k}. \end{equation}
\end{thm}

{\bf Proof.} Let
$$ \tilde H = \begin{pmatrix} & & \Gamma_2 & \Gamma_1 \cr & & \ddots & \Gamma_2 \cr & & &
\cr & & & \end{pmatrix} . $$
Then it follows from \eqref{YW} that
$$ \begin{pmatrix} I & \tilde H \cr \tilde{H}^* & I \end{pmatrix} \begin{pmatrix}
B \cr \tilde{D} \end{pmatrix} = - \begin{pmatrix} \Gamma \cr 0 \end{pmatrix}, $$
where $$ \Gamma =\begin{pmatrix} \Gamma_0 \cr \Gamma_1 \cr \vdots \end{pmatrix} , \tilde D= \begin{pmatrix} \vdots \cr
D_{-2} \cr D_{-1} \end{pmatrix} . $$
But then it follows that \eqref{gj} is equivalent to the equation
\begin{equation}\label{gj1} \Gamma_j = \begin{pmatrix} 0 & Z_{j+1} \end{pmatrix} \begin{pmatrix} I & \tilde H \cr
\tilde{H}^* & I \end{pmatrix}^{-1} \begin{pmatrix} \Gamma \cr 0 \end{pmatrix} , j \le -1 ,
\end{equation}
where
$$ Z_k = \begin{pmatrix} \cdots & \Gamma_{k+1} & \Gamma_k \end{pmatrix} . $$
But this coincides exactly with the iterative process described in \cite{AAK4} (see also \cite[Section 2.2]{Peller}), and thus
the conclusion follows from there.

For the alternative construction of $\Gamma_j$, use that \eqref{YW2} implies that
$$ \begin{pmatrix} I & \tilde H \cr \tilde{H}^* & I \end{pmatrix} \begin{pmatrix}
\tilde{A} \cr C \end{pmatrix} = - \begin{pmatrix} 0 \cr \hat{\Gamma} \end{pmatrix}, $$
where $$ \hat{\Gamma} =\begin{pmatrix} \vdots \cr \Gamma_1^* \cr \Gamma_0^* \end{pmatrix} ,
\tilde A= \begin{pmatrix}
A_1 \cr A_2 \cr \vdots \end{pmatrix} . $$
But then \eqref{gj2} is equivalent to the equality
$$ \Gamma_j^* = \begin{pmatrix} \hat{Z}_{j+1} & 0 \end{pmatrix}
\begin{pmatrix} I & \tilde H \cr \tilde{H}^* & I \end{pmatrix}^{-1}
\begin{pmatrix} 0 \cr \hat{\Gamma} \end{pmatrix} , $$
with
$$ \hat{Z}_k = \begin{pmatrix} \Gamma_{k} & \Gamma_{k+1} & \cdots \end{pmatrix} , $$
which yields the same sequence of operators $\Gamma_k$, $k \le -1$, as in \eqref{gj1}.
$\Box$

We now come to the main result in this section.

\begin{thm}\label{main2}
Let $\gamma_{ij} \in L({\mathcal H}, {\mathcal K})$, $i,j\ge 0$, be given so that the little
Hankel operator $h_\gamma : l^2_{\mathcal K} ({-{\mathbb N}_0 \times -{\mathbb N}_0}) \to
l^2_{\mathcal H} ({{\mathbb N}_0 \times {\mathbb N}_0})$ defined via
$$ h_\gamma = \begin{pmatrix} & & \Gamma_1 & \Gamma_0 \cr & & \ddots & \Gamma_1 \cr & & &
\cr & & & \end{pmatrix} \ , \
\Gamma_j = \begin{pmatrix} & & \gamma_{j1} & \gamma_{j0} \cr & & \ddots & \gamma_{j1} \cr & & &
\cr & & & \end{pmatrix} , $$
is a strict contraction. Put
$$ \Phi = P_{{\mathbb N}_0 \times {\mathbb N}} \oplus P_{-{\mathbb N} \times -{\mathbb N}_0}
\begin{pmatrix} I & h_\gamma \cr h_\gamma^* & I \end{pmatrix}
P^*_{{\mathbb N}_0 \times {\mathbb N}} \oplus P^*_{-{\mathbb N} \times -{\mathbb N}_0}
= P_{{\mathbb N} \times {\mathbb N}_0} \oplus P_{-{\mathbb N}_0 \times -{\mathbb N}}
\begin{pmatrix} I & h_\gamma \cr h_\gamma^* & I \end{pmatrix}
P^*_{{\mathbb N} \times {\mathbb N}_0} \oplus P^*_{-{\mathbb N}_0 \times -{\mathbb N}}, $$
$$
\Phi_1 = P_{{\mathbb N}_0 \times {\mathbb N}_0} \oplus P_{-{\mathbb N} \times -{\mathbb N}}
\begin{pmatrix} I & h_\gamma \cr h_\gamma^* & I \end{pmatrix}
P^*_{{\mathbb N} \times {\mathbb N}_0} \oplus P^*_{-{\mathbb N}_0 \times -{\mathbb N}}
= P_{{\mathbb N}_0 \times {\mathbb N}} \oplus P_{-{\mathbb N} \times -{\mathbb N}_0}
\begin{pmatrix} I & h_\gamma \cr h_\gamma^* & I \end{pmatrix}
P^*_{{\mathbb N} \times {\mathbb N}} \oplus P^*_{-{\mathbb N}_0 \times -{\mathbb N}_0}, $$
$$
\Phi_2 = P_{{\mathbb N}_0 \times {\mathbb N}_0} \oplus P_{-{\mathbb N} \times -{\mathbb N}}
\begin{pmatrix} I & h_\gamma \cr h_\gamma^* & I \end{pmatrix}
P^*_{{\mathbb N}_0 \times {\mathbb N}} \oplus P^*_{-{\mathbb N} \times -{\mathbb N}_0}
= P_{{\mathbb N} \times {\mathbb N}_0} \oplus P_{-{\mathbb N}_0 \times -{\mathbb N}}
\begin{pmatrix} I & h_\gamma \cr h_\gamma^* & I \end{pmatrix}
P^*_{{\mathbb N} \times {\mathbb N}} \oplus P^*_{-{\mathbb N}_0 \times -{\mathbb N}_0}, $$
where the projection
$P_K : l^2(M) \to l^2(K), \ K \subseteq M,$ is defined by $P_k((\eta_j)_{j\in M})
= (\eta_j)_{j\in K} . $
Suppose that
\begin{equation}\label{comm}
\Phi_1 \Phi^{-1} \Phi_2^* = \Phi_2^* \Phi^{-1} \Phi_1.
\end{equation}
Then there exist $\gamma_{ij} \in L ({\mathcal H} )$, $(i,j) \in ({\mathbb Z} \times -{\mathbb N})
\cup (-{\mathbb N} \times {\mathbb Z})$, so that the operator matrix
$$ (\gamma_{i-j,k-l})_{i,j,k,l \in {\mathbb Z}}: l^2_{\mathcal K} ({\mathbb Z} \times
{\mathbb Z}) \to l^2_{\mathcal H} ({\mathbb Z} \times
{\mathbb Z}) $$
is a strict contraction. Equivalently, the essentially bounded function $f \sim \sum_{i,j\in
{\mathbb Z}} \gamma_{ij}z^iw^j$ satisfies $\| f \|_\infty < 1 $.
\end{thm}

{\bf Proof.} We start by applying Theorem \ref{Nehari} to construct $\Gamma_j$, $j \le -1$, via
\eqref{gj} or, equivalently, \eqref{gj2}, yielding the strict contraction
$$ (\Gamma_{i-j})_{i,j\in{\mathbb Z}} : l^2_{l^2_{\mathcal K} (-{\mathbb N}_0 )} ({\mathbb Z})
\to l^2_{l^2_{\mathcal H} ({\mathbb N}_0 )} ({\mathbb Z}). $$
The main step in the proof is to show that \eqref{comm} implies that $\Gamma_j$, $j \le -1$, are also Hankel; that is, they are of the form
$$ \Gamma_j = \begin{pmatrix} & & \gamma_{j1} & \gamma_{j0} \cr & & \ddots & \gamma_{j1} \cr & & &
\cr & & & \end{pmatrix}, j \le -1, $$
for some operators $\gamma_{ij}$, $j\ge 0$, $i \le -1$. To show this we need to prove the following claim.

{\bf Claim.} \it Equation \eqref{comm} implies that the operators $D_j$, $j \le -1$,
in \eqref{YW} are of the form
$$ D_j = \begin{pmatrix} & \vdots & \vdots & \vdots \cr \cdots & * & * & * \cr \cdots & * & * & *
\cr \cdots & 0 & 0 & * \end{pmatrix} : l^2_{\mathcal K } (-{\mathbb N}_0 ) \to
l^2_{\mathcal K } (-{\mathbb N}_0 ). $$
Similarly, \eqref{comm} implies that $A_j$ in \eqref{YW2} is of the form
$$ A_j = \begin{pmatrix} * & 0 & 0 & \cdots \cr * & * & * & \cdots
\cr  * & * & * & \cdots \cr \vdots & \vdots & \vdots & \end{pmatrix} :
l^2_{\mathcal H } ({\mathbb N}_0 ) \to
l^2_{\mathcal H } ({\mathbb N}_0 ). $$
\rm

{\bf Proof of Claim.} It is not hard to see that $\Phi_i \Phi^{-1}$ and $\Phi^{-1} \Phi_i$, $i=1,2$,
have a certain companion type form (variations of the ones in the proof of Theorem \ref{main}).
For instance,
$$ \Phi_1 \Phi^{-1} = \begin{pmatrix} \hat{S} & Q \cr 0 & S \end{pmatrix},
\Phi^{-1} \Phi_1 = \begin{pmatrix} Z & \hat{Q} \cr 0 & \hat{Z} \end{pmatrix}, $$
where $\hat{S}$ and $\hat{Z}$ have an infinite companion form
$$ \hat{S} = \begin{pmatrix} * & *&  \cdots \cr I & & \cr & I & \cr & & \ddots \end{pmatrix},
\hat{Z} = \begin{pmatrix} \ddots & & & \vdots \cr & I & & * \cr & & I & * \end{pmatrix}, $$
the operators $S$ and $Z$ are shifts
$$ S = \begin{pmatrix} \ddots & & & \vdots \cr & I & & 0 \cr & & I & 0 \end{pmatrix},
Z = \begin{pmatrix} 0 & 0&  \cdots \cr I & & \cr & I & \cr & & \ddots \end{pmatrix},$$
and $Q$ and $\hat{Q}$ are zero except for the first block row and last
block column, respectively:
$$ Q =  \begin{pmatrix} * & *&  \cdots \cr 0 & 0 & \cdots \cr \vdots & \vdots & \end{pmatrix},
\hat{Q} =  \begin{pmatrix} & \vdots & \vdots \cr \cdots  &0 & * \cr \cdots & 0 & * \end{pmatrix}. $$
But then, viewing $R:= \Phi_2^* \Phi^{-1} \Phi_1 = \Phi_1 \Phi^{-1} \Phi_2^* $ in the four possible ways $ (\Phi_2^* \Phi^{-1}) \Phi_1,$ $   \Phi_2^* (\Phi^{-1} \Phi_1), $ $( \Phi_1 \Phi^{-1}) \Phi_2^*,$ $  \Phi_1 (\Phi^{-1} \Phi_2^*)$ one easily deduces that
$$ \Phi_1 \Phi^{-1}\Phi_2^* = \Phi_2^* \Phi^{-1} \Phi_1 = P_{{\mathbb N}_0 \times {\mathbb N}} \oplus P_{-{\mathbb N} \times
-{\mathbb N}_0} \begin{pmatrix} I & h_\gamma \cr h_\gamma^* & I \end{pmatrix}
P^*_{{\mathbb N} \times {\mathbb N}_0} \oplus P^*_{-{\mathbb N}_0 \times -{\mathbb N}} . $$
Multiplying the above equation on the left with $0 \oplus P_{-{\mathbb N} \times \{ 0 \} }$
and on the right with $0 \oplus P^*_{\{ 0 \} \times -{\mathbb N}}$ gives that
$ Y W^{-1} U = X $, where $U,W,X$ and $Y$ are defined via
$$ Y =0 \oplus P_{-{\mathbb N} \times \{ 0 \} } \begin{pmatrix} I & h_\gamma \cr h_\gamma^* & I \end{pmatrix} P^*_{{\mathbb N}_0 \times {\mathbb N}_0} \oplus P^*_{-{\mathbb N} \times -{\mathbb N}}, $$
$$ W =  P_{{\mathbb N}_0 \times {\mathbb N}_0} \oplus P_{-{\mathbb N} \times -{\mathbb N}}
\begin{pmatrix} I & h_\gamma \cr h_\gamma^* & I \end{pmatrix}
P^*_{{\mathbb N}_0 \times {\mathbb N}_0} \oplus P^*_{-{\mathbb N} \times -{\mathbb N}}, $$
$$ U = P_{{\mathbb N}_0 \times {\mathbb N}_0} \oplus P_{-{\mathbb N} \times -{\mathbb N}}
\begin{pmatrix} I & h_\gamma \cr h_\gamma^* & I \end{pmatrix} 0 \oplus P^*_{\{ 0 \} \times -{\mathbb N}} , $$
and
$$ X = 0 \oplus P_{-{\mathbb N} \times \{ 0 \} } \begin{pmatrix} I & h_\gamma \cr h_\gamma^* & I \end{pmatrix} 0 \oplus P^*_{\{ 0 \} \times -{\mathbb N}} . $$
View the operator
$$ M = \begin{pmatrix} I & h_\gamma P^*_{(-{\mathbb N}_0 \times -{\mathbb N}_0) \setminus
\{ (0,0) \} } \cr P_{(-{\mathbb N}_0 \times -{\mathbb N}_0) \setminus
\{ (0,0) \} } h_\gamma^* & I \end{pmatrix} $$
after permutation as the operator matrix
$$ \begin{pmatrix} * & Y & X \cr * & W & U \cr * &* & * \end{pmatrix} $$
acting on
$$ [0 \oplus l^2(-{\mathbb N} \times \{ 0 \} )] \oplus
[l^2({\mathbb N}_0 \times {\mathbb N}_0 ) \oplus l^2 (-{\mathbb N} \times -{\mathbb N}) ]
\oplus [0 \oplus l^2(\{ 0 \} \times -{\mathbb N}) ] . $$
Then the equality $ Y W^{-1} U = X $ together with Proposition \ref{3x3} gives that
$$ (0 \oplus P_{-{\mathbb N} \times \{ 0 \} }) M^{-1} ( 0 \oplus P^*_{\{ 0 \} \times -{\mathbb N}} ) =0 .$$
This exactly yields the required zeros in $D_j$, $j \le -1$.

The proof of the zeros in $A_j$, $j \ge 1$, is similar. This proves the claim. $\Box$

Following the claim, we may now write $D_j$ and $A_j$ as
\begin{equation}\label{djaj}
D_j = \begin{pmatrix} \widetilde{D}_j & q_j \cr 0 & \delta_j \end{pmatrix} , j \le -1;
\ \ A_j =  \begin{pmatrix} \alpha_j & 0 \cr r_j & \widehat{A}_j \end{pmatrix} , j \ge 1.
\end{equation}
Write
$$ \Gamma_j = \begin{pmatrix} \widetilde{\Gamma}_j & \begin{matrix} \gamma_{j0} \cr
\gamma_{j1} \cr \vdots \end{matrix} \end{pmatrix} ,  \Gamma_j = \begin{pmatrix}
\begin{matrix} \cdots & \gamma_{j1} & \gamma_{j0}  \end{matrix} \cr \widehat{\Gamma}_j \end{pmatrix}
. $$
Note that $\widetilde{\Gamma}_j = \widehat{\Gamma}_j$, $j \ge 0$. Observe that due
to \eqref{djaj}, equation \eqref{YW} implies
$$ \begin{pmatrix} I & \widetilde{h}_\gamma \cr \widetilde{h}_\gamma^* & I \end{pmatrix}
\begin{pmatrix} B \cr \widetilde{D} \end{pmatrix} = \begin{pmatrix} 0 \cr \widetilde{\Delta}
\end{pmatrix} , $$
with $$  \widetilde{h}_\gamma = (\widetilde{\Gamma}_{i-j} )_{i\in{\mathbb N}_0, j\in -{\mathbb N}_0},
\widetilde{D} = \begin{pmatrix} \vdots \cr \widetilde{D}_{-2} \cr \widetilde{D}_{-1} \cr I
\end{pmatrix} , \widetilde{\Delta} = \begin{pmatrix} \vdots \cr 0 \cr 0 \cr \widetilde{\Delta}_0
\end{pmatrix}, $$
where $\widetilde{\Delta}_0$ is obtained from $\Delta_0$ by removing the last row
and column; that is $ \Delta_0 = \begin{pmatrix}  \widetilde{\Delta}_0 & * \cr * & * \end{pmatrix} $.
Moreover, if we define
$$ \widetilde{\Gamma}_j = - \sum_{k=1}^\infty \widetilde{\Gamma}_{j+k} \widetilde{D}_{-k} , j
\le -1 , $$
then we have that $\widetilde{\Gamma}_j$ corresponds to $\Gamma_j$ without the last column for
$j\le -1$ as well. In other words,
$$ \Gamma_j = \begin{pmatrix} \widetilde{\Gamma}_j & * \end{pmatrix} , j \le -1 . $$
Likewise, due to the form of $A_j$, we have that $\widehat{A}_j$ may be constructed from
\eqref{YW2} with $\Gamma_j$ replaced by $\widehat{\Gamma}_j$. Moreover, if we define
$$ \widehat{\Gamma}_j^* = - \sum_{k=1}^\infty \widehat{\Gamma}_{j+k}^* \widehat{A}_j, j \le -1, $$
then we have that $\Gamma_j = \begin{pmatrix} * \cr \widehat{\Gamma}_j \end{pmatrix}$, $ j \le
-1 $. But since $\widehat{\Gamma}_j = \widetilde{\Gamma}_j$, $j \ge 0$, we obtain from
Theorem \ref{Nehari} that
$$ \widetilde{\Gamma}_j = - \sum_{k=1}^\infty \widetilde{\Gamma}_{j+k} \widetilde{D}_{-k}
=  - \sum_{k=1}^\infty \widehat{\Gamma}_{j+k} \widetilde{D}_{-k} =
- \sum_{k=1}^\infty \widehat{A}_j^* \widehat{\Gamma}_{j+k} = \widehat{\Gamma}_j, \
j \le -1 . $$
Since
$$ \Gamma_j =  \begin{pmatrix} * \cr \widehat{\Gamma}_j \end{pmatrix} = \begin{pmatrix} \widetilde{\Gamma}_j & * \end{pmatrix} , j \le -1 , $$
it now follows that $\Gamma_j$, $j \le -1$, is Hankel.

The last step in the proof is to recognize that
$$ \| (\Gamma_{i-j})_{i,j\in {\mathbb Z}} \| < 1 $$
implies that the Hankel
$ (H_{i-j})_{i \in {\mathbb N}_0, j \in -{\mathbb N}_0} $ is a strict contraction,
where
$$ H_i = (\gamma_{p-q, i })_{p,q \in {\mathbb Z}} , i \ge 0 . $$
But now it follows that $H_i = (\gamma_{p-q, i })_{p,q \in {\mathbb Z}}$, $i \le -1$, exist
so that $(H_{i-j})_{i,j= -\infty}^\infty $ is a strict contraction. $\Box$

\bibliographystyle{plain} 

\begin{thebibliography}{10}

\bibitem{AAK4}
V.~M. Adamjan, D.~Z. Arov, and M.~G. Kre\u{\i}n.
\newblock Infinite {H}ankel block matrices and related problems of extension.
\newblock {\em Izv. Akad. Nauk Armjan. SSR Ser. Mat.}, 6(2-3):87--112, 1971.

\bibitem{ACC}
Gr. Arsene, Zoia Ceau{\c{s}}escu, and T.~Constantinescu.
\newblock Schur analysis of some completion problems.
\newblock {\em Linear Algebra Appl.}, 109:1--35, 1988.

\bibitem{CP}
A.~Calderon and R.~Pepinsky.
\newblock On the phases of fourier coefficients for positive real periodic
  functions.
\newblock {\em Computing Methods and the Phase Problem in $X$-Ray Crystal
  Analysis (R. Pepinsky, ed.)}, pages 339--346, 1950.

\bibitem{CS2}
Mischa Cotlar and Cora Sadosky.
\newblock Two distinguished subspaces of product {BMO} and {N}ehari-{AAK}
  theory for {H}ankel operators on the torus.
\newblock {\em Integral Equations Operator Theory}, 26(3):273--304, 1996.

\bibitem{DGK2}
Philippe Delsarte, Yves~V. Genin, and Yves~G. Kamp.
\newblock Planar least squares inverse polynomials. {I}. {A}lgebraic
  properties.
\newblock {\em IEEE Trans. Circuits and Systems}, 26(1):59--66, 1979.

\bibitem{DGK3}
Philippe Delsarte, Yves~V. Genin, and Yves~G. Kamp.
\newblock Half-plane {T}oeplitz systems.
\newblock {\em IEEE Trans. Inform. Theory}, 26(4):465--474, 1980.

\bibitem{DGK4}
Philippe Delsarte, Yves~V. Genin, and Yves~G. Kamp.
\newblock{A simple proof of {R}udin's multivariable stability theorem},
\newblock{\em IEEE Trans. Acoust. Speech Signal Process}, 28(6):701--705, (1980), 

\bibitem{DGK5}
Philippe Delsarte, Yves~V. Genin, and Yves~G. Kamp.
\newblock Half-plane minimization of matrix-valued quadratic functionals.
\newblock {\em SIAM J. Algebraic Discrete Methods}, 2(2):192--211, 1981.

\bibitem{Di}
Bradley~W. Dickinson.
\newblock Two-dimensional markov spectrum estimates need not exist.
\newblock {\em IEEE Trans. Inform. Theory}, 26:120--121, 1980.

\bibitem{DG79a}
Harry Dym and Israel Gohberg.
\newblock Extensions of matrix valued functions with rational polynomial
  inverses.
\newblock {\em Integral Equations Operator Theory}, 2(4):503--528, 1979.

\bibitem{DG79}
Harry Dym and Israel Gohberg.
\newblock Extensions of band matrices with band inverses.
\newblock {\em Linear Algebra Appl.}, 36:1--24, 1981.

\bibitem{DG5}
Harry Dym and Israel Gohberg.
\newblock Extensions of kernels of {F}redholm operators.
\newblock {\em J. Analyse Math.}, 42:51--97, 1982/83.

\bibitem{FL}
Sarah~H. Ferguson and Michael~T. Lacey.
\newblock A characterization of product {BMO} by commutators.
\newblock {\em Acta Math.}, 189(2):143--160, 2002.

\bibitem{FS}
Sarah~H. Ferguson and Cora Sadosky.
\newblock Characterizations of bounded mean oscillation on the polydisk in
  terms of {H}ankel operators and {C}arleson measures.
\newblock {\em J. Anal. Math.}, 81:239--267, 2000.

\bibitem{FF}
Ciprian Foias and Arthur~E. Frazho.
\newblock {\em The commutant lifting approach to interpolation problems}.
\newblock Birkh\"auser Verlag, Basel, 1990.

\bibitem{FF90}
Ciprian Foias and Arthur~E. Frazho.
\newblock {\em The commutant lifting approach to interpolation problems},
  volume~44 of {\em Operator Theory: Advances and Applications}.
\newblock Birkh\"auser Verlag, Basel, 1990.

\bibitem{GK1}
Y.~Genin and Y.~Kamp.
\newblock Counterexample in the least-squres inverse stabilizatin of {\rm 2d}
  recursive filters.
\newblock {\em Electron Lett.}, 11:330--331, 1975.

\bibitem{GW}
Jeffrey~S. Geronimo and Hugo~J. Woerdeman.
\newblock Positive extensions, {F}ej\'er-{R}iesz factorization and
  autoregressive filters in two variables.
\newblock {\em Ann. of Math. (2)}.
\newblock to appear.

\bibitem{GH}
I.~Gohberg and G.~Heinig.
\newblock Inversion of finite {T}oeplitz matrices consisting of elements of a
  noncommutative algebra.
\newblock {\em Rev. Roumaine Math. Pures Appl.}, 19:623--663, 1974.

\bibitem{GKW1}
I.~Gohberg, M.~A. Kaashoek, and H.~J. Woerdeman.
\newblock The band method for positive and contractive extension problems.
\newblock {\em J. Operator Theory}, 22:109--155, 1989.

\bibitem{GKW2}
I.~Gohberg, M.~A. Kaashoek, and H.~J. Woerdeman.
\newblock The band method for positive and contractive extension problems: {A}n
  alternative version and new applications.
\newblock {\em Integral Equations and Operator Theory}, 12:343--382, 1989.

\bibitem{GGK2}
Israel Gohberg, Seymour Goldberg, and Marinus~A. Kaashoek.
\newblock {\em Classes of linear operators. {V}ol. {II}}, volume~63 of {\em
  Operator Theory: Advances and Applications}.
\newblock Birkh\"auser Verlag, Basel, 1993.

\bibitem{HE}
Henry Helson.
\newblock {\em Lectures on invariant subspaces}.
\newblock Academic Press, New York, 1964.

\bibitem{HL}
Henry Helson and David Lowdenslager.
\newblock Prediction theory and {F}ourier series in several variables.
\newblock {\em Acta Math.}, 99:165--202, 1958.

\bibitem{HL2}
Henry Helson and David Lowdenslager.
\newblock Prediction theory and {F}ourier series in several variables. {I}{I}.
\newblock {\em Acta Math.}, 106:175--213, 1961.

\bibitem{L-APK}
Hanoch Lev-Ari, Sydney~R. Parker, and Thomas Kailath.
\newblock Multidimensional maximum-entropy covariance extension.
\newblock {\em IEEE Trans. Inform. Theory}, 35(3):497--508, 1989.

\bibitem{LM}
J.~S. Lim and N.~A. Malik.
\newblock A new algorithm for two-dimensional maximum entropy power spectrum
  estimation.
\newblock {\em IEEE Trans. Acoust., Speech, Signal Processing}, 29:401--412,
  1981.

\bibitem{Page}
Lavon~B. Page.
\newblock Applications of the {S}z.-{N}agy and {F}oia\c s lifting theorem.
\newblock {\em Indiana Univ. Math. J.}, 20:135--145, 1970/1971.

\bibitem{Peller}
Vladimir~V. Peller.
\newblock {\em Hankel operators and their applications}.
\newblock Springer Monographs in Mathematics. Springer-Verlag, New York, 2003.

\bibitem{RO}
R.~Rochberg.
\newblock Toeplitz operators on weighted ${H}^p$ spaces.
\newblock {\em Indiana Univ. Math. J.}, 26:291--298, 1977.

\bibitem{Rud}
Walter Rudin.
\newblock The extension problem for positive-definite functions.
\newblock {\em Illinois J. Math.}, 7:532--539, 1963.

\bibitem{Wo2}
H.~J. Woerdeman.
\newblock {\em Matrix and operator extensions}.
\newblock Stichting Mathematisch Centrum Centrum voor Wiskunde en Informatica,
  Amsterdam, 1989.

\end{thebibliography}
\def\polhk#1{\setbox0=\hbox{#1}{\ooalign{\hidewidth
  \lower1.5ex\hbox{`}\hidewidth\crcr\unhbox0}}} \def\cprime{$'$}
  \def\cprime{$'$}

\end{document}